\documentclass[a4paper,11pt]{article}
\usepackage{amssymb}
\usepackage{amsmath}
\usepackage{amsfonts}
\usepackage[mathscr]{eucal}
\usepackage{enumerate}
\usepackage{amscd}
\usepackage{indentfirst}
\usepackage{enumerate}
\usepackage{txfonts, textcomp}
\newtheorem{thm}{Theorem}[section]
\newtheorem{prop}[thm]{Proposition}
\newtheorem{cor}[thm]{Corollary}
\newtheorem{dfn}[thm]{Definition}

\newtheorem{ex}[thm]{Example}
\newtheorem{rmk}{Remark}[section]
\numberwithin{equation}{section}
\allowdisplaybreaks[4]

\title{EQUIVALENCE OF THE CATEGORIES OF MODULES OVER LIE ALGEBROIDS}
\author{Yuji HIROTA\thanks{hirota@rs.tus.ac.jp}
\\Department of Mathematics\\Tokyo University of Science}
\date{}
\begin{document}
\maketitle 

\begin{abstract}
We study geometric representation theory of Lie algebroids. A new equivalence relation 
for integrable Lie algebroids is introduced and investigated. 
It is shown that two equivalent 
Lie algebroids have equivalent categories of infinitesimal actions of Lie algebroids. 
As an application, it is also shown that the Hamiltonian categories for 
gauge equivalent Dirac structures are equivalent as categories. 
\end{abstract}

\noindent {\bf Mathematics Subject Classification~(2000).} ~53D17 

\section{Introduction}
Poisson geometry is considered to be intermediate between differential geometry 
and noncommutative geometry in the sense that it provides us with powerful techniques 
to study many geometric objects related to noncommutative algebras. 
%In several research domain of Poisson geometry, the theory of geometric Morita equivalence 

If $(Q,\,\Pi_Q)$ and $(P,\,\Pi_P)$ are Poisson manifolds, then 
a Poisson map $J:Q\to P$ induces a Lie algebra homomorphism by 
\begin{equation}\label{sec1:eqn}
C^\infty(P)\longrightarrow \mathfrak{X}(Q)\subset \mathrm{End}\,\bigl(C^\infty(Q)\bigr), 
\quad f\longmapsto -\Pi_Q(\,\cdot,\,J^*df\,).
\end{equation}
From (\ref{sec1:eqn}), $C^\infty(Q)$ can be regarded as a $C^\infty(P)$-module. 
This observation enables oneself to study geometric objects by connecting with 
a theory in algebra like Morita equivalence (refer to H. Bursztyn and A. Weinstein \cite{BWpoi05}
for further discussion). Geometric 
Morita equivalence, which is introduced by P. Xu~\cite{Xmor91}, 
plays a central role in Poisson geometry as Morita equivalence of $C^\star$-algebras does 
in noncommutative geometry. One of the remarkable properties is that Morita equivalence implies 
the equivalence of the categories of modules over Poisson manifolds: 
for integrable Poisson manifold $P$, the category of modules over $P$ is the category whose 
objects are complete symplectic realizations of $P$ and whose morphisms are symplectic maps 
between complete symplectic realizations commuting with the realizations. 
This is just the analogy with Morita equivalence in algebra, first studied by 
K. Morita~\cite{Mdua58}. 

As is well-known, Poisson maps are always associated with Lie algebroid actions of cotangent 
bundles: 
\begin{equation}\label{sec1:eqn2}
\varGamma^\infty(T^*P)\longrightarrow\mathfrak{X}(Q),\quad \alpha\longmapsto 
-\Pi_Q(\,\cdot,\,J^*\alpha\,).
\end{equation}
The Lie algebra homomorphism (\ref{sec1:eqn2}) can be considered to be the representation of 
$\varGamma^\infty(T^*P)$ on $C^\infty(Q)$. More generally, if $A\to M$ is a Lie algebroid, 
the infinitesimal action of $A$ on a smooth map $f:N\to M$ induces the representation of $\varGamma^\infty(A)$ 
on $C^\infty(N)$: 
\begin{equation}\label{sec1:eqn3}
\varGamma^\infty(A)\longrightarrow\mathfrak{X}(N)\subset \mathrm{End}\bigl(C^\infty(N)\bigr). 
\end{equation}
Here, a natural question arises: what is an equivalence relation between Lie algebroids 
which implies an equivalence of the categories associating with Lie algebroid actions? 
%such as Morita equivalence for integrable Poisson manifolds?

In this paper, we give a solution to the above question, that is, 
we introduce an equivalence relation for integrable Lie algebroids, 
called strong Morita equivalence, and show that the category 
consisting of the infinitesimal actions of  
Lie algebroids is invariant under strong Morita equivalent. 
Furthermore, applying the result to Dirac geometry, 
we partially recover the well-known proposition in H. Bursztyn and M. Crainic \cite{BCdir08}. 
This study gives a general description of Morita equivalence for Poisson manifolds 
from the viewpoint of Lie algebroid, 
and is expected to have a connection with the study of quasi-Hamiltonian symmetry 
through the question presented by A. Weinstein \cite{Wgeo05}.

The paper is organized as follows: in Section 2, we review the basics of Lie algebroids, 
including Lie algebroid morphisms and the construction of Lie algebroid from a given Lie groupoid. 
Section 3 is devoted to the study of the infinitesimal actions of Lie algebroids. 
The new equivalence relation for integrable Lie algebroids 
is introduced and discussed. 
In Section 4, we show that 
the category of the infinitesimal actions of Lie algebroid is invariant under 
strong Morita equivalence, and show also that two gauge equivalent Dirac structures are 
strong Morita equivalent. Lastly, we find that 
the Hamiltonian categories for gauge equivalent Dirac structures are equivalent each other, 
by using the main theorem. 

Throughout the paper, the set of smooth sections of a smooth vector bundle $E\to X$ is denoted 
by $\varGamma^\infty(X)$. Especially, we write $\mathfrak{X}(M)$ for $\varGamma^\infty(TX)$ 
when $E=TX$. The space of smooth functions on a smooth manifold $M$ is denoted by $C^\infty(M)$. 
\vspace{0.5cm}

\noindent{\bf Acknowledgments.}~The author would like to thank Ping Xu for 
helpful advice on the study and Akira Yoshioka for various support. 
He also wishes to thank Keio University for the hospitality while part of the work was being done. 
%%%%%%%%%%%%%%%%%%%%%%%%%%%%%%%%%%%%%%%%%%%%%%%%%%%%%%%%%%%%%%%%%%%%%%%%%%%%%%%%%%
\section{Basic terminologies of Lie algebroids}
\subsection{Lie algebroids}

Let $M$ be a smooth manifold. A Lie algebroid over $M$ is a smooth vector bundle $A\to M$ 
with a bundle map $\rho:E\to TM$, called the anchor map, and 
a Lie bracket $\llbracket\cdot,\,\cdot\rrbracket$ on the space $\varGamma^\infty(A)$ of 
smooth sections of $A$ such that 
\begin{equation}\label{sec2:eqn1}
 \llbracket\alpha,\,f\beta\rrbracket 
            = \bigl(\rho(\alpha)f\bigr)\,\beta \,+\, f\llbracket\alpha,\,\beta\rrbracket
\end{equation}
for any $f\in C^\infty(M)$ and $\alpha,\,\beta\in \varGamma^\infty(A)$. 
We denote a Lie algebroid by the triple $(A\to M,\,\llbracket\cdot,\,\cdot\rrbracket,\,\rho)$ or, 
simply by $A$, and use the notation $A^{-}$ for a Lie algebroid $A$ with the opposite bracket. 
\vspace{0.4cm}

The anchor map of a Lie algebroid $A$ is a Lie algebra homomorphism. Indeed, from 
(\ref{sec2:eqn1})~ and the Jacobi identity, it follows that 
\begin{align*}
 0 &= \llbracket\llbracket\alpha,\,\beta\rrbracket,\, f\gamma\rrbracket \,+\,
      \llbracket\llbracket\beta,\,f\gamma\rrbracket,\, \alpha\rrbracket \,+\, 
      \llbracket\llbracket f\gamma,\,\alpha\rrbracket,\, \beta\rrbracket\\
   &= f\llbracket\llbracket\alpha,\,\beta\rrbracket,\, \gamma\rrbracket \,+\, 
   \left(\rho\bigl(\llbracket\alpha,\,\beta\rrbracket\bigr)f\right)\gamma\\ 
   &\quad \,+\, 
   f\llbracket\llbracket\beta,\,\gamma\rrbracket,\, \alpha\rrbracket \,+\, 
   \bigl(\rho(\beta)f\bigr)\llbracket\gamma,\,\alpha\rrbracket 
      \,-\, \bigl(\rho(\alpha)f\bigr)\llbracket\beta,\,\gamma\rrbracket 
           \,-\, \left(\rho(\alpha)\bigl(\rho(\beta)f\bigr)\right)\gamma\\
  &\quad \,+\, f\llbracket\llbracket\gamma,\,\alpha\rrbracket,\, \beta\rrbracket 
       \,-\, \bigl(\rho(\beta)f\bigr)\llbracket\gamma,\,\alpha\rrbracket 
         \,+\, \bigl(\rho(\alpha)f\bigr)\llbracket\beta,\,\gamma\rrbracket
           \,+\, \left(\rho(\beta)\bigl(\rho(\alpha)f\bigr)\right)\gamma\\
  &= \left(\bigl(\rho\bigl(\llbracket\alpha,\,\beta\rrbracket\bigr) - 
   \llbracket\rho(\alpha),\,\rho(\beta)\rrbracket\bigr)f\right)\gamma
\end{align*}
for any $f\in C^\infty(M)$ and $\alpha,\,\beta,\,\gamma\in \varGamma^\infty(A)$. 
Therefore, we have 
$\rho\bigl(\llbracket\alpha,\,\beta\rrbracket\bigr) 
  =\llbracket\rho(\alpha),\,\rho(\beta)\rrbracket$.
\begin{ex}
A Lie algebra is a Lie algebroid over a point. 
\end{ex}
\begin{ex}(Tangent algebroids)~
A tangent bundle $TM$ of a smooth manifold $M$ is a Lie algebroid over $M$$\mathrm{:}$ 
the anchor map is 
the identity map $\mathrm{id}_{TM}$, 
and the Lie bracket is the usual Lie bracket of vector fields. This Lie algebroid 
is called a tangent algebroid. 
\end{ex}
\begin{ex}(Cotangent algebroids)~
If $(P,\,\Pi)$ is a Poisson manifold, then a cotangent bundle $T^*P$ is a Lie algebroid$\mathrm{:}$ 
the anchor map is the map $\Pi^\sharp$ induced from $\Pi$, 
\[\Pi^\sharp: T^*P\longrightarrow TP,\quad \alpha\longmapsto 
\bigl\{\, \beta\mapsto \langle\beta,\,\Pi^\sharp(\alpha)\rangle = \Pi(\beta,\,\alpha)\, \bigr\}\]
and the Lie bracket is given by 
\[\llbracket\alpha,\,\beta\rrbracket = \mathcal{L}_{\Pi^\sharp(\alpha)}\beta
 - \mathcal{L}_{\Pi^\sharp(\beta)}\alpha + d\bigl(\Pi(\alpha,\,\beta)\bigr),\] 
where $\mathcal{L}_{\Pi^\sharp(\alpha)}\beta$ stands for the Lie derivative on $\beta$ 
along $\Pi^\sharp(\alpha)$. 
The Lie algebroid $(T^*P\to P,\,\llbracket\cdot,\,\cdot\rrbracket,\,\Pi^\sharp)$ is called a 
cotangent algebroid. 
\end{ex}
\begin{ex}(Transformation algebroids)~
Given an action $\varrho:\mathfrak{g}\to \mathfrak{X}(M)$ of a Lie algebra 
$(\mathfrak{g},\,[\cdot,\,\cdot])$ 
on a smooth manifold $M$, one can associate to it the Lie algebroid structure$\mathrm{:}$ 
the vector bundle is the trivial bundle $M\times\mathfrak{g}\to M$, the anchor map 
$\rho$ is given by 
$\rho\, (p,\,V)\to \left(\varrho(V)\right)_p\in T_pM,\,(\forall p\in M,\,V\in \mathfrak{g})$ and 
the Lie bracket on $\varGamma^\infty(M\times\mathfrak{g}) = C^\infty(M,\,\mathfrak{g})$ is defined as 
\[\llbracket U,\,V\rrbracket (p) := [U(x),\,V(x)] + \bigl(\varrho(U(x))\bigr)_p(V) 
- \bigl(\varrho(V(x))\bigr)_p(U).\]
This Lie algebroid is called a transformation algebroid, and denoted by $\mathfrak{g}\ltimes M$,
for short. 
\end{ex}
\begin{ex}\label{sec2:example of Lie algebroid}(Dirac structures)~
Let us consider a vector bundle $TM\oplus T^*M$ over a smooth manifold $M$. 
We endow the vector bundle with a bilinear operation 
\[\langle\cdot,\,\cdot\rangle :
\varGamma^\infty(TM\oplus T^*M)\times\varGamma^\infty(TM\oplus T^*M)\to C^\infty(M) 
\]
defined as 
\[
\langle\,(X,\alpha),\,(Y,\beta)\,\rangle := \beta(U) + \alpha(V), 
\]
and a skew-symmetric bracket
\[\llbracket\cdot,\,\cdot\rrbracket : 
\varGamma^\infty(TM\oplus T^*M)\times\varGamma^\infty(TM\oplus T^*M)\to
\varGamma^\infty(TM\oplus T^*M)
\]
defined as 
\[
\llbracket(U,\alpha),\,(V,\beta)\rrbracket := (\,[U,V],\,\mathcal{L}_U\beta - \mathrm{i}_Vd\alpha\,). 
\]
A subbundle $D_M\subset TM\oplus T^*M$ is called a Dirac structure if $D_M$ satisfies the following
three conditions:
\begin{enumerate}[\rm\quad(1)]
 \item $\left.\langle\cdot,\,\cdot\rangle\right|_{D_M} \equiv 0$; 
 \item $D_M$ has rank equal to $\dim(M)$;
 \item 
  $\llbracket\varGamma^\infty(D_M),\,\varGamma^\infty(D_M)\rrbracket\subset \varGamma^\infty(D_M)$.
\end{enumerate}
We call 
a pair $(M, D_M)$ of a smooth manifold $M$ and a Dirac structure $D_M\subset TM\oplus T^*M$ 
a Dirac manifold. 
A Dirac structure $D_M$, with the restriction 
of Courant bracket and the anchor map, is verified easily to be a Lie algebroid. 
We refer to H. Bursztyn and M. Crainic~\cite{BCdir05}, and \cite{BWpoi05} 
for further discussions of Dirac structures and Courant algebroids. 
\end{ex}

%%%%%%%%%%%%%%%%%%%%%%%%%%%%%%%%%%%%%%%%%%%%%%%%%%%%%%%%%%%%%%%%%%%%%%%%%%%%%%%%%
\subsection{Lie algebroid morphisms and the pull-back Lie algebroids}
%%%%%%%%%%%%%%%%%%%%%%%%%%%%%%%%%%%%%%%%%%%%%%%%%%%%%%%%%%%%%%%%%%%%%%%%%%%%%%%%

Let $(A_1\to M_1,\,\llbracket\cdot,\,\cdot\rrbracket^1,\,\rho_1)$ and 
$(A_2\to M_2,\,\llbracket\cdot,\,\cdot\rrbracket^2,\,\rho_2)$ be Lie algebroids. 
A Lie algebroid morphism from $A_1$ to $A_2$ is a vector bundle morphism 
$\Phi:A_1\to A_2$ such that 
\begin{equation}\label{Lie morphism 1}
   \rho_2\bigl(\Phi(\alpha)\bigr) = \Phi_*\bigl(\rho_1(\alpha)\bigr),\quad 
  \bigl(\forall \alpha\in \varGamma^\infty(A_1)\bigr),
\end{equation}
and, for any smooth sections $\alpha,\,\beta\in \varGamma^\infty(A_1)$ written in the forms 
\begin{equation}\label{Lie morphism 2}
\Phi(\alpha) = \sum_i\xi_i(\gamma_i\circ \varphi),\quad 
\Phi(\beta) = \sum_j\eta_j(\delta_j\circ \varphi),
\end{equation}
where $\xi_i,\,\eta_j\in C^\infty(M_1)$ and $\gamma_i,\,\delta_j\in \varGamma^\infty(A_2)$, 
\begin{equation}\label{Lie morphism 3}
   \Phi\bigl(\llbracket\alpha,\,\beta\rrbracket^1\bigr) = 
    \sum_{i,j}\xi_i \eta_j (\llbracket\gamma_i,\,\delta_j\rrbracket^2\circ\Phi) 
    + \sum_j\bigl(\mathcal{L}_{\rho_1(\alpha)}\eta_j\bigr) (\delta_j\circ\Phi) 
    - \sum_i\bigl(\mathcal{L}_{\rho_1(\beta)}\xi_i\bigr) (\gamma_i\circ\Phi) 
\end{equation}
are satisfied (see K. Mackenzie \cite{Mgen05}). 
Here, we denote the base map of $\Phi$ 
by $\underline{\Phi}$. 

\begin{prop}
If a vector bundle morphism $\Phi:A_1\to A_2$ is the Lie algebroid morphism, 
then there exists a subbundle 
\[R\subset \left.(A_1\times A_2)\right |_{\mathrm{Gr}(\Phi)}\]
which satisfies the following conditions{\rm :} 
 \begin{enumerate}[\quad\rm (1)]
  \item For any $z\in \mathrm{Gr}(\Phi)$, 
   $(\rho_1\times \rho_2)\bigl((R)_z\bigr)\subset T_z\bigl(\mathrm{Gr}(\Phi)\bigr)\mathrm{;}$
  \item For any $\alpha,\,\beta\in \varGamma^\infty(A_1\times A_2)$ such that 
 $\left.\alpha\right|_{\mathrm{Gr}(\Phi)},\,\left.\beta\right|_{\mathrm{Gr}(\Phi)}
 \in \varGamma^\infty(R)$, $\left.\llbracket\,\alpha, \,\beta\,\rrbracket\right|_{\mathrm{Gr}(\Phi)} 
 \in \varGamma^\infty(R)$,
 \end{enumerate}
where $\llbracket\cdot,\,\cdot\rrbracket 
= \bigl(\,\llbracket\cdot,\,\cdot\rrbracket^1,\,\llbracket\cdot,\,\cdot\rrbracket^2\,\bigr)$.
\end{prop}
({\it Proof})~Suppose that $\varphi:A_1\to A_2$ is a Lie algebroid morphism. 
Define the vector bundle $R\subset \left.(A_1\times A_2)\right |_{\mathrm{Gr}(\Phi)}$ as 
\[R \,=\, \coprod_{p\in M_1}\Bigl\{\,(a,\,\Phi(a))\,\bigm|\, 
a\in (A_1)_p\,\Bigr\}.
\]
Using (\ref{Lie morphism 1}), we have 
\begin{align*}
(\rho_1\times\rho_2)\bigl(a,\,\Phi(a)\bigr) 
 &= \Bigl(\,\rho_1(a),\,\rho_2\bigl(\Phi(a)\bigr)\,\Bigr)\\ 
 &= \Bigl(\,\rho_1(a),\,\Phi_*\bigl(\rho_1(a)\bigr)\,\Bigr) 
 \in T_p\bigl(\mathrm{Gr}(\Phi)\bigr).
\end{align*}
That is, the condition (1) holds.

For $\alpha,\,\beta\in\varGamma^\infty(A_1)$ which we assume to satisfy (\ref{Lie morphism 2}), 
we define the smooth sections $\widehat{\alpha},\,\widehat{\beta}$ 
of $A_1\times A_2\to M_1\times M_2$ as  
\[\widehat{\alpha}_{(p,\Phi(p))} 
:= \bigl(\alpha_p,\,\Phi(\alpha_p)\bigr) \in (R)_{(p,\Phi(p))},\quad
 \widehat{\beta}_{(p,\Phi(p))} := \bigl(\beta_p,\,\Phi(\beta_p)\bigr) 
 \in (R)_{(p,\varphi(p))}. 
\]
From (\ref{sec2:eqn1}) and (\ref{Lie morphism 3}), it follows that 
\begin{equation*}
 {\llbracket\,\Phi(\alpha),\,\Phi(\beta)\,\rrbracket^2}_{\Phi(p)} 
 = \Phi\bigl({\llbracket\,\alpha,\,\beta\,\rrbracket^1}_p\bigr).
\end{equation*}
This leads us to the condition (2). 
\qquad\qquad\qquad\qquad\qquad\quad\qquad\qquad\qquad\qquad\quad\qquad\qquad\quad$\Box$
\vspace*{0.4cm}

The Lie algebroid morphism $\Phi:A_1\to A_2$ is said to 
be a Lie algebroid isomorphism if $\Phi$ is an isomorphism of vector bundles. 
If there exists the Lie algebroid isomorphism from $A_1$ to $A_2$, we write $A_1\cong A_2$. 
\vspace*{0.5cm}

Let $(A\to M,\,\llbracket\cdot,\cdot\rrbracket,\,\rho)$ be a Lie algebroid and 
$f:M'\to M$ a smooth map from a smooth manifold $M'$ to $M$. 
Assume that the differential of $f$ is transversal to the anchor map 
$\rho:A\to TM$ in the sense that 
\[\mathrm{Im}\,\rho_{f(x)} \,+\, \mathrm{Im}\,(df)_x \,=\, 
         T_{f(x)}M,\quad (\forall x\in M').\]
Here, $\mathrm{Im}\,\rho_{f(x)}$ stands for the image of $\rho_{f(x)}$.

This assumption leads us to the following condition:
 \begin{equation}\label{sec1:transversality}
 \mathrm{Im}\,(I_x\times\rho_{f(x)}) \,+\, T_{(x.f(x))}
  \bigl(\mathrm{Gr}(f)\bigr)
  \,=\, T_xM'\oplus T_{f(x)}M,\quad (\forall x\in M'), 
 \end{equation}
where $I_x$ means the identity map on $T_xM'$.
The condition (\ref{sec1:transversality})~ensures that the preimage 
\begin{equation}\label{pullback}
(I\times\rho)^{-1}T\bigl(\mathrm{Gr}(f)\bigr) \,=\, 
\coprod_{x\in M'}\Bigl\{\,(V,\,\alpha)\,\bigm|\, V\in T_xM',\,\alpha\in A_{f(x)},\,
(df)_x(V) = \rho(\alpha)\,\Bigr\}
\end{equation} 
is a smooth subbundle of $(\left.TM'\times A)\right|_{\mathrm{Gr}(f)}$. 
The vector bundle (\ref{pullback}) over $\mathrm{Gr}(f)\cong M'$ has the structure of Lie algebroid 
whose anchor map is the natural projection $\mathrm{proj}_1$. 
This vector bundle is called a pull-back of Lie algebroid and denoted by $f^!A$ 
(see P. Higgins and K. Mackenzie \cite{HMalg90}). 
\vspace*{0.5cm}

Let $\Phi_1:A_1\to A$ and $\Phi_2:A_2\to A$ be Lie algebroid morphisms. We denote each base map 
by $\underline{\Phi}_1:M_1\to M$ and $\underline{\Phi}_2:M_2\to M$. 
Suppose that the following conditions: 
 \begin{enumerate}[\quad\rm (1)]
  \item $\mathrm{Im}\,(\Phi_1)_p \,+\,\mathrm{Im}\,(\Phi_2)_q = A_r,\quad 
    \bigl(\,r=\Phi_1(p)=\Phi_2(q)\,\bigr)\mathrm{;}$
  \item The map $\underline{\Phi}_1\times\underline{\Phi}_2$ is transversal to the submanifold 
   $\bigtriangleup = \{\,(m,\,m)\,|\,m\in M\,\}\subset M\times M$:
    \[\mathrm{Im}\,\bigl((d\underline{\Phi}_1)_p\times (d\underline{\Phi}_2)_q\bigr) \,+\, 
           T_{(r,r)}\bigtriangleup \,=\, T_{(r,r)}(M\times M)\]
 \end{enumerate}
are satisfied. Then, one can obtain the Lie algebroid 
\[A_1\times_A A_2 := \coprod_{(p,q)\in M_1\times_M M_2}\Bigl\{\,(a,\,b)\,\bigm|\, 
a\in (A_1)_p,\,b\in (A_2)_q,\,\Phi_1(a) = \Phi_2(b) \,\Bigr\}\] 
over $M_1\times_M M_2 = \bigl\{(p,\,q)\in M_1\times M_2 \,|\, \underline{\Phi}_1(p) 
= \underline{\Phi}_2(q)\bigr\}$, 
whose Lie bracket $\llbracket\cdot,\,\cdot\rrbracket$ is given by $\llbracket\cdot,\,\cdot\rrbracket 
:= \bigl(\,\llbracket\cdot,\,\cdot\rrbracket^1,\,\llbracket\cdot,\,\cdot\rrbracket^2\,\bigr)$, 
and whose anchor map $\widehat{\rho}:A_1\times_A A_2\to T(M_1\times_M M_2)$ is defined as 
$\widehat{\rho}(a,\,b) := (\rho_1(a),\,\rho_2(b))$. 
We call this Lie algebroid the fibered product. 
The pull-back of a Lie algebroid $f^!A$ discussed can be 
the fibered product of two Lie algebroid morphisms $f_*:TM'\to TM$ and $\rho:A\to TM$. 
Hence, a fibered product Lie algebroid is a pull-back Lie algebroid in a general sense. 
%%%%%%%%%%%%%%%%%%%%%%%%%%%%%%%%%%%%%%%

\subsection{The Lie algebroid of a Lie groupoid}

Let $\varGamma\rightrightarrows M$ be a Lie groupoid with an identity section 
$\boldsymbol{\varepsilon}$, a source map $\boldsymbol{s}$ and a target map $\boldsymbol{t}$. 
Denote by $\mathcal{A}\,(\varGamma)\to M$ the vector bundle consisting of tangent spaces 
to $\boldsymbol{s}$-fibers at $X$:
\[ 
\left.\mathcal{A}\,(\varGamma)\right|_p = \ker(d\boldsymbol{s})_{\boldsymbol{\varepsilon}(p)}
\quad (p\in M). 
\] 
For any $\gamma\in\varGamma$, the differential of the right translation $R_{\gamma}$ 
by $\gamma$ induces a map 
\[
 (dR_\gamma)_{\boldsymbol{\varepsilon}(\gamma')} : 
 T_{\boldsymbol{\varepsilon}(\gamma')}
  \left(\boldsymbol{s}^{-1}\bigl(\boldsymbol{t}(\gamma)\bigr)\right)\longrightarrow 
 T_{\boldsymbol{\varepsilon}(\gamma')}\left(\boldsymbol{s}^{-1}\bigl(\boldsymbol{s}(\gamma)\bigr)\right), 
\]
where $\gamma'=\boldsymbol{t}(\gamma)$. 
By the map, any smooth section $\alpha\in\varGamma^\infty(\mathcal{A}\,(\varGamma))$  
gives rise to a right-invariant vector field  
\begin{equation}\label{sec2:right-invariant}
 \widehat{\alpha}_\gamma := (dR_\gamma)_{\boldsymbol{\varepsilon}(\gamma')}
   (\alpha_{\boldsymbol{\varepsilon}(\gamma')})\quad (\gamma\in \varGamma) 
\end{equation}
on $\varGamma$ (see \cite{Mgen05}). 
Therefore, $\varGamma^\infty(\mathcal{A}\,(\varGamma))$ inherits the Lie bracket from 
$\mathfrak{X}(\varGamma)$. 
One verifies that the vector bundle $\mathcal{A}\,(\varGamma)\to M$ with the above 
Lie bracket and the bundle map 
$d\boldsymbol{t}:\mathcal{A}\,(\varGamma)\to TM$ becomes a Lie algebroid. 
A Lie algebroid $A\to M$ is said to be integrable if there exists a Lie groupoid 
$\varGamma\rightrightarrows M$ whose Lie algebroid $\mathcal{A}(\varGamma)\to M$ 
is isomorphic to $A$ as Lie algebroid. 

%Before proceeding with the next section, we give a short review of Lie's 
%third theorem for Lie algebroids, here. 
%For further details, we refer to M. Crainic and R. J. Fernandes~\cite{CFint03}, \cite{CFint04}. 
%\begin{dfn}
% Let $A\overset{\pi}{\to} M$ be a Lie algebroid. 
% An $A$-path is a smooth path $a:I\to A$ which projects 
% to a base path $\pi\circ a:I\to M$ such that 
% \[
%  \rho\bigl(a(t)\bigr) \,=\, \frac{d}{dt}\pi\bigl(a(t)\bigr)\quad (\forall t\in I). 
% \]
% Here, $I=[0,\,1]$ stands for the unit interval. 
%\end{dfn}

%%%%%%%%%%%%%%%%%%%%%%%%%%%%%%%%%%%%%%%%%%%%%%%%%%%%%%%%%%%%%%%%%%%%%%%%%%%%%%%%%%

\section{Infinitesimal actions of Lie algebroids and strong Morita equivalence}

%%%%%%%%%%%%%%%%%%%%%%%%%%%%%%%%%%%%%%%%%%%%%%%%%%%%%%%%%%%%%%%%%%%%%%%%%%%%

We begin this section by recalling the actions of Lie algebroids. 
A Lie algebroid right (left) action of $(A\to M,\,\llbracket\cdot, \cdot\rrbracket,\,\rho)$ 
on a smooth manifold $N$ consists of a map $\mu:N\to M$ called the 
momentum map and a Lie algebra (anti-) homomorphism 
$\xi:\varGamma^\infty(A)\to \mathfrak{X}(N)$ which satisfy 
\begin{equation}\label{sec3:eqn1}
\rho\bigl(\alpha_{\mu(q)}\bigr) = (d\mu)_q\bigl(\xi(\alpha)\bigr) \quad (\forall q\in N)
\end{equation}
for any $\alpha\in\varGamma^\infty(A)$, and 
\begin{equation}\label{sec3:eqn2}
\rho (f\alpha) = (\mu^*f)\, \xi(\alpha) \quad (\forall f\in C^\infty(M)).
\end{equation}
The right action of $A$ is alternatively called the infinitesimal action of $A$. 
The action is said to be complete if $\xi(\alpha)$ is a complete vector field whenever 
$\alpha\in \varGamma^\infty(A)$ has compact support. 

%The action of $A$ is said to be free if for each $x\in X$ the mapping 
%$\varGamma^\infty(\mu^*A)\ni\alpha\mapsto \xi(\alpha)_x\in T_xX$ is injective. 
%Then, the distribution spanned by $\xi(\alpha)$ is of constant rank. 
%On the other hand, 
%the action is said to be transitive if for each $x\in X$ the mapping 
%$\varGamma^\infty(\mu^*A)\ni\alpha\mapsto \xi(\alpha)_x\in T_xX$ is surjective onto $T_xX$. 

\begin{ex}
Let $\mathfrak{g}$ be a Lie algebra. A Lie algebra action of $\mathfrak{g}$ on $M$ 
is thought of a Lie algebroid action of $\mathfrak{g}\to \{\ast\}$ on $M\to\{\ast\}$. 
\end{ex}

\begin{ex}\label{sec3:Poisson map}
Any Poisson map $J:Q\to P$ and a cotangent algebroid $T^*P$ over $P$ is a Lie algebroid action 
by (\ref{sec1:eqn2}). 
\end{ex}

\begin{ex}
Any smooth manifold $X$ is thought of a Lie algebroid action 
of a trivial Lie algebroid $\{\ast\}\to\{\ast\}$ on a map $X\to \{\ast\}$. 
We call this action a trivial action. 
\end{ex}

\begin{ex}\label{sec3:example4}
Given a Lie algebroid $A\to M$ with a surjective submersion $J:X\to M$ which satisfy 
\begin{equation}\label{sec4:eqn3}
({J}^!E)_{(x, J(x))} \,\cap\, T_xX = \{\boldsymbol{0}\}\quad (\forall x\in X), 
\end{equation}
we have the right action of Lie algebroid $\varGamma^\infty(J^*A)\to \mathfrak{X}(X)$ by 
$\alpha\mapsto u$, where $u\in T_xX$ is the element such that $(u,\alpha)\in({J}^!E)_{(x, J(x))}$. 
We remark that the element $u$ is uniquely determined by {\rm (\ref{sec4:eqn3})}. Indeed, 
if $\alpha$ and $u, u'$ are the elements such that $(u,\alpha)\in({J}^!E)_{(x, J(x))}$ and 
$(u',\alpha)\in({J}^!E)_{(x, J(x))}$, then $(u-u', \boldsymbol{0})$. 
So, we have $u=u'$. 
\end{ex}

\begin{ex}\label{sec4:ex2}
Let us assume that a Lie algebroid $A\to M$ is integrable and $\varGamma\rightrightarrows M$ 
be the Lie groupoid integrating $A$. As noted in Section 2, the fiber of $A$ over $x\in M$ 
is the subspace $\ker\,(d\boldsymbol{s})_{\varepsilon(x)}$ of $T_{\varepsilon(x)}\varGamma$, and 
the anchor is given by $d\boldsymbol{t}:A\subset T\varGamma\to TM$. 
Given any section $\alpha\in\varGamma^\infty(A)$, the formula (\ref{sec2:right-invariant})
defines a right invariant vector field. 
The map $\xi$ which assigns the right invariant vector field $\widehat{\alpha}$ on $\varGamma$ 
to $\alpha\in\varGamma^\infty(A)$ is shown to be a Lie algebra homomorphism and 
satisfy (\ref{sec3:eqn1}) and (\ref{sec3:eqn2}). Therefore, the map 
$\varGamma^\infty(A)\to \mathfrak{X}(\varGamma)$ defines a right action of $A$ on 
$\boldsymbol{t}:\varGamma\to M$. Similarly to this case, one can obtain a left action of $A$ 
on $\boldsymbol{t}:\varGamma\to M$ by defining as 
$\varGamma^\infty(A)\ni\alpha\mapsto -\widehat{\alpha}\in\mathfrak{X}(\varGamma)$. 
\end{ex}

\begin{prop}\label{sec4:prop}
Let $(A\to M, \llbracket\cdot,\cdot\rrbracket, \rho)$ be a Lie algebroid and 
$J:X\to M$ a smooth map. Suppose that $J$ is a surjective submersion. 
Then, we have a Lie algebroid action of $A$ on $X/\mathcal{F}$, 
where $X/\mathcal{F}$ is the space of leaves induced from $J$.
\end{prop}

\noindent({\it Proof})~
Since $J$ is a surjective submersion, the space $X$ has a foliation $\mathcal{F}$ whose leaves are $J$-fibers. 
We consider the space of leaves $X/\mathcal{F}$ and  
a map $\overline{J}: X/\mathcal{F}\to M$ given by $\overline{J}(\overline{x})=J(x)~(\forall x\in X)$. 
For any $\alpha_{J(x)}\in A_{J(x)}~(x\in X)$, there exists $u_x\in T_xX$ such that 
$(dJ)_x(u_x)=\rho(\alpha_{J(x)})$. 
A vector field $u=\{u_x\}_{x\in X}\in\mathfrak{X}(X)$ is $J$-related to 
$\rho(\alpha)\in\mathfrak{X}(M)$: $dJ\circ u = \rho(\alpha)\circ J$. 
We define a map $\xi:\varGamma^\infty(A)\to \mathfrak{X}(X/\mathcal{F})$ as 
\begin{equation*}
A_{J(x)}\longrightarrow T_{\overline{x}}(X/\mathcal{F}), \quad 
\alpha_{J(x)}\longmapsto \overline{u}_{\overline{x}}:=(d\pi)_x(u_x),
\end{equation*} 
where $\pi$ stands for a natural projection $\pi:X\to X/\mathcal{F},\, x\to \overline{x}$. 
Let $\xi(\alpha)=\overline{u}$ and $\xi(\beta)=\overline{v}$ for 
$\alpha, \beta\in\varGamma^\infty(A)$. The vector fields $\overline{u}$ and 
$\overline{v}$ on $X/\mathcal{F}$ 
are $\pi$-related to $u$ and $v$, respectively. It follows from this that 
$[\xi(\alpha), \xi(\beta)] = \overline{[u, v]}$. 
On the other hand, we take a vector field $w$ on $X$ such that 
$\rho\,(\llbracket\alpha, \beta\rrbracket)\circ J=dJ\circ w$. 
Since the anchor map $\rho$ is a Lie algebra homomorphism (see Section 2), we have 
\begin{align*}
 w(J^*g) &= (dJ\circ w)f = \bigl([\rho(\alpha),\,\rho(\beta)]\circ J\bigr)g 
  = \bigl([dJ\circ u,\,dJ\circ v]\bigl)g\\
 &= \bigl(dJ\circ [u,\,v])\bigr)g = [u,\,v](J^*g) 
\end{align*}
for any $g\in C^\infty(M)$. In other words, it holds that $w=[u, v]$ on each $J$-fiber. 
Hence, we have $\xi(\llbracket\alpha, \beta\rrbracket) = \overline{[u, v]}$. 
These result in that the map $\xi$ is a Lie algebra homomorphism. 
It is shown easily that $\xi$ also satisfies (\ref{sec3:eqn1}) and (\ref{sec3:eqn2}). 
\qquad\qquad\qquad\qquad\qquad\qquad\qquad\qquad\qquad\qquad\qquad\qquad $\Box$ 

\begin{rmk}\label{sec4:rmk}
If a Lie algebroid $A$ acts on $\mu:N\to M$, then a pull-back vector bundle $\mu^*A\to N$ 
has a Lie algebroid structure whose anchor is the action map. We refer to \cite{Mgen05} for 
further details. 
\end{rmk}
\vspace{0.5cm}

From the definition of the Lie algebroid action, the space $C^\infty(N)$ can be regarded as a 
$\varGamma^\infty(A)$-module. In other words, one can think 
the actions of Lie algebroids of the modules over Lie algebroids. 
We define a right (left) module over a Lie algebroid $A$ to be the right (resp. left) action of $A$ 
whose momentum map is a surjective submersion. 
A right (left) module over $A$ is said to be complete if the right (resp. left) action is 
complete. 

\begin{ex}
The action of $T^*P\to P$ given by 
$\varGamma^\infty(T^*P)\ni\alpha\mapsto \Pi_Q(\,\cdot,\,J^*\alpha\,)\in\mathfrak{X}(Q)$ 
is a left module over $T^*P$ (see {\rm (\ref{sec1:eqn2})}). 
\end{ex}

\begin{ex}
The Lie algebroid action of $A$ in Proposition $\ref{sec4:prop}$ is the right module over $A$. 
\end{ex}

\begin{ex}\label{sec3:groupoid action}
Let $\varGamma_1\rightrightarrows \varGamma_0$ be a Lie groupoid. 
Let us take points $x\in \varGamma_0$ and $h\in\varGamma_1$ such that 
$\boldsymbol{t}(h)=x$. 
For any smooth section $\alpha\in\varGamma^\infty(\mathcal{A}(\varGamma_1))$, 
we consider a smooth curve $\gamma$ in $\boldsymbol{s}^{-1}(x)$ which satisfies 
\[
\left.\frac{d}{dt}\right|_{t=0}\gamma = \widehat{\alpha}_x,\quad \text{\rm and} \quad 
\gamma(0) = \varepsilon(x). 
\]
Since $\boldsymbol{s}(\gamma(t)) = x =\boldsymbol{t}(h)$ for each $t\in\mathbb{R}$, 
a smooth curve $t\mapsto \gamma(t)\cdot h$ can be defined. 
Then, the map 
\begin{equation}\label{sec3:groupoid action1}
\varGamma^\infty(\mathcal{A}(\varGamma_1))\longrightarrow \mathfrak{X}(\varGamma_1),\quad 
\alpha\longmapsto -\left\{\left.\frac{d}{dt}\right|_{t=0}\gamma(t)\cdot h\right\}_{h\in \varGamma_1} 
\end{equation}
defines a left module $\boldsymbol{t}:\varGamma_1\to \varGamma_0$ 
over $\mathcal{A}(\varGamma_1)\to \varGamma_0$.  

On the other hand, let us 
consider a smooth curve $\delta$ in $\boldsymbol{t}^{-1}(x)$ which satisfies 
\[
\left.\frac{d}{dt}\right|_{t=0}\delta = \widehat{\beta}_x,\quad \text{\rm and} \quad 
\delta(0) = \varepsilon(x). 
\] 
for $x\in \varGamma_0$ and $g\in\varGamma_1$ such that 
$\boldsymbol{s}(g)=x$, and for 
any smooth section $\beta\in\varGamma^\infty(\mathcal{A}(\varGamma_1))$. 
Then, the map defined as 
\begin{equation}\label{sec3:groupoid action2}
\varGamma^\infty(\mathcal{A}(\varGamma_1))\longrightarrow \mathfrak{X}(\varGamma_1),\quad 
\beta\longmapsto \left\{\left.\frac{d}{dt}\right|_{t=0}g\cdot\delta(t)\right\}_{h\in \varGamma_1} 
\end{equation}
is a right module $\boldsymbol{s}:\varGamma_1\to \varGamma_0$ 
over $\mathcal{A}(\varGamma_1)\to \varGamma_0$.  
\end{ex}
\vspace*{0.5cm}

Suppose that we are given a right $A$-module $J:X\to M$ by 
\begin{equation*}
 \xi: \varGamma^\infty(A)\longrightarrow \mathfrak{X}(X), \quad 
      (A)_{J(x)}\ni \alpha_{J(x)}\longmapsto \xi\,(\alpha_{J(x)})\in T_xX 
\end{equation*}
and a left $A$-module $K:Y\to M$ by 
\begin{equation*}
 \eta: \varGamma^\infty(A)\longrightarrow \mathfrak{X}(Y), \quad 
       (A)_{K(y)}\ni \alpha_{K(y)}\longmapsto \eta\,(\alpha_{K(y)}) \in T_yY. 
\end{equation*}
Take the fiber product 
\[
X\times_{M_2}Y = \Bigl\{\,(x, y)\in X\times Y\,\bigm|\,J(x)=K(y)\,\Bigr\}, 
\]
then, we have a right $A$-action on $\widehat{J}:X\times_{M}Y\to M$ by 
\begin{equation}\label{sec4:eqn4}
 (A)_{J(x)}\ni \alpha_{J(x)}\longmapsto 
\bigl(\xi\,(\alpha)_x,\,-\eta\,(\alpha)_y\bigr)\in T_{(x,y)}\bigl(X\times_{M}Y\bigr). 
\end{equation} 
This action gives rise to a singular distribution $D=\{D_{(x,y)}\}$ on $X\times_{M}Y$
\[X\times_{M}Y\ni (x, y)\longmapsto 
D_{(x,y)} := 
\Bigl\{\,\bigl(\xi_2(\alpha)_x,\,-\eta_2(\alpha)_y\bigr)
\,\bigm|\, \alpha\in \varGamma^\infty(A_2)\,\Bigr\}\subset T_{(x,y)}\bigl(X\times_{M_2}Y\bigr).
\]
The distribution $D$ turns out to be integrable since 
the action (\ref{sec4:eqn4}) is thought of the anchor map of the fibered product 
$J^*A\times_{A}K^*A\to X\times_{M}Y$~
(see Remark \ref{sec4:rmk}, and 8.1.4 in J.-P. Dufour and N. T. Zung~\cite{DZpoi05}). 
We denote by $X\otimes_{A} Y$ 
the space of leaves $(X\times_{M}Y)/\mathcal{A}$ obtained from $D$. 
\vspace*{0.5cm}

\begin{dfn}\label{sec3:dfn} 
 Two Lie algebroids $A_1\to M_1$ and $A_2\to M_2$ are said to be quasi-equivalent 
 if there exists 
 a smooth manifold $X$ together with surjective submersions 
 $J_k:X\to M_k~(k=1,2)$ and a pair $(L_1, L_2)$ of 
 subbundles $L_1$ of ${J_1}^!A_1$ and $L_2$ of ${J_2}^!A_2^{-}$ which satisfy the 
 following conditions{\rm :}
  \begin{enumerate}[\quad\rm(A)]
   \item $(L_k)_{(x, J_k(x))} \,\cap\, (T_xX\oplus\{\boldsymbol{0}\}) = \{\boldsymbol{0}\}$\quad 
         {\rm for any} $x\in X$\quad $(k=1,2){\rm ;}$ 
  %\item $(L_k)_{(x, J_k(x))} \,\cap\, (\{\boldsymbol{0}\}\oplus E_{J_k(x)}) = 
         %\{\boldsymbol{0}\}$\quad 
         %{\rm for any} $x\in X$\quad $(k=1,2){\rm ;}$ 
   \item $\mathrm{pr}_1\Bigl((L_1)_{(x,J_1(x))}\Bigr) = 
                T_x\Bigl(J_2^{-1}\bigl(J_2(x)\bigr)\Bigr)$ \quad
         {\rm and} \quad 
         $\mathrm{pr}_1\Bigl((L_2)_{(x, J_2(x))}\Bigr) = 
                T_x\left(J_1^{-1}\bigl(J_1(x)\bigr)\right)$\quad 
         $(\forall x\in X){\rm ;}$ 
   \item $\mathrm{pr}_2\Bigl((L_1)_{(x,J_1(x))}\Bigr) = (A_1)_{J_1(x)}$ \quad
         {\rm and} \quad $\mathrm{pr}_2\Bigl((L_2)_{(x,J_2(x))}\Bigr) = (A_2)_{J_2(x)}$\quad 
         $(\forall x\in X){\rm;}$
  \end{enumerate}
 where $\mathrm{pr}_1$ and $\mathrm{pr}_2$ are natural projections from $TX\times A_i~(i=1,2)$ 
 to the first component $TX$ and the second component $E_i$, respectively. 
\end{dfn}

\begin{ex}\label{sec3:example11}
Suppose that integrable Poisson manifolds $P_1$ and $P_2$ are Morita equivalent in the 
sense of Xu~\cite{Xmor91} each other, 
that is, there exists a symplectic manifold $S$ together with two surjective submersions 
$P_1\overset{\tau_1}{\leftarrow}S\overset{\tau_2}{\rightarrow}P_2$ 
such that 
 \begin{enumerate}[\quad\rm(1)]
  \item $\tau_1$ is a complete Poisson map and 
         $\tau_2$ is a complete anti-Poisson map{\rm ;}
  \item each $\tau_k$ has connected, simply connected fibers\quad $(k=1,2)${\rm ;} 
  \item $\ker(d\tau_1)_z=\bigl(\ker(d\tau_2)_z\bigr)^\perp$\quad {\rm and}\quad 
        $\ker(d\tau_2)_z=\bigl(\ker(d\tau_1)_z\bigr)^\perp$ \quad $(\forall z\in S)$.
 \end{enumerate}
Then, the cotangent algebroids $T^*P_1\to P_1$ and $T^*P_2\to P_2$ 
are quasi-equivalent: to verify this claim, let us take subbundles 
$L_1\subset {\tau_1}^!(T^*P_1)$ and $L_2\subset {\tau_2}^!(T^*P_2)^-$ as 
\[
 L_1 = \bigl\{\,(\Pi_S^\sharp({\tau_1}^*\alpha),\,\alpha)\,|\,\alpha\in T^*P_1\,\bigr\}
\]
and 
\[
 L_2 = \bigl\{\,(\Pi_S^\sharp({\tau_2}^*\beta),\,\beta)\,|\,\beta\in (T^*P_2)^-\,\bigr\},
\]
respectively, where $\Pi^\sharp_S$ 
stands for the bundle map induced by the symplectic Poisson structure 
$\Pi_S\in \varGamma^\infty(\wedge^2TS)$.  
The condition (A) and (C) in Definition \ref{sec3:dfn} are easily checked. 
The condition (B) follows from (3) in the above 
that fibers of $\tau_1,\,\tau_2$ are symplectically orthogonal to one another. 
\end{ex}

The statement similar to Example \ref{sec3:example11} holds in a more general setting. 
Let $D_{M_1}$ and $D_{M_2}$ be Dirac structures over $M_1$ and $M_2$, respectively. 
A smooth map $F:M_1\to M_2$ is called a forward Dirac map if it holds that 
\[
 (D_{M_2})_{F(m)} = \left\{\,\bigl((dF)_mU,\, \beta\bigr)\in T_mM_1\oplus T_{F(m)}M_2\,\bigm|
        \,\bigl(U,\,(dF)_m^*\beta\bigr)\in (D_{M_1})_m\,\right\}
\]
for any point $m\in M_1$. In addition, a forward Dirac map $F:(M_1, D_{M_1})\to (M_2, D_{M_2})$ 
is called a strong Dirac map if 
\begin{equation}\label{sec3:strong Dirac}
\ker(dF)_m \cap \ker(D_{M_1})_m 
= \{\boldsymbol{0}\}\quad (\forall m\in M_1)
\end{equation}
is satisfied, 
where $\ker(D_{M_1})_m = (D_{M_1})_m\cap T_mM_1$ (see H. Bursztyn and M. Crainic ~\cite{BCdir08}). 

\begin{rmk}
 A strong Dirac map is alternatively called a Dirac realization in \cite{BCdir05}. 
\end{rmk}

A strong Dirac map $F:(M_1, D_{M_1})\to (M_2, D_{M_2})$ induces a map 
\begin{equation}\label{sec3:Dirac action}
 \zeta : \varGamma^\infty(D_{M_2})\longrightarrow \mathfrak{X}(M_1),\quad 
  (V,\, \beta) \longmapsto \hat{V}, 
\end{equation}
where $\hat{V}$ is a tangent vector such that $V=F_*\hat{V}$ 
which is determined uniquely by the condition 
(\ref{sec3:strong Dirac}). 
The map $\zeta$ defines an infinitesimal actions of the Lie algebroid $D_{M_2}$ 
(see Proposition 2.3 in \cite{BCdir08}). 
A strong Dirac map $F$ is said to be complete if the infinitesimal action $\zeta$ is 
complete. 
As noted in Example \ref{sec2:example of Lie algebroid}, 
Dirac structures are regarded as Lie algebroids. The following proposition states 
the sufficient condition for two Dirac structures to be quasi-equivalent. 

\begin{prop}
Two Dirac structures $D_{M_1}$ and $D_{M_2}$ are quasi-equivalent if 
there exists a Dirac manifold $(N, D_N)$ together with surjective submersions 
$(M_1, D_{M_1})\overset{F_1}{\leftarrow}N\overset{F_2}{\rightarrow}(M_2, D_{M_2}^-)$ satisfying
 \begin{enumerate}[\quad\rm(1)]
  \item each $F_k$ is a strong Dirac map\quad $(k=1,2)${\rm ;}
  \item $\mathrm{pr}_1\bigl((\varLambda_1)_n\bigr) = \ker\,(dF_2)_n$\quad {\rm and}\quad 
   $\mathrm{pr}_2\bigl((\varLambda_2)_n\bigr) = \ker\,(dF_1)_n$\quad $(\forall n\in N)$,
 \end{enumerate} 
 where $(\varLambda_k)_n:=(D_N)_n\cap \bigl(T_nN\oplus \mathrm{Im}\,(dF_k)_n^*\bigr)\quad (k=1,2)$. 
\end{prop}

\noindent({\it Proof})~ We define subbundles $L_1\subset {F_1}^!D_{M_1}$ and 
$L_2\subset {F_2}^!D_{M_2}^-$ over $N$ as 
\[
 L_1 := \coprod_{n\in N}\Bigl\{\,\bigl(u;\, (dF_1)_n(u),\, \beta\bigr)\, 
 \bigm|\,u\in T_nN,\, \beta\in T_{F_1(n)}M_1,\, \bigl(u,\, (df_1)_n^*(\beta)\bigr)\in (D_N)_n\,\Bigr\}
\]
and 
\[
 L_2 := \coprod_{n\in N}\Bigl\{\,\bigl(u;\, (dF_2)_n(u),\, \beta\bigr)\, 
 \bigm|\,u\in T_nN,\, \beta\in T_{F_2(n)}M_1,\, \bigl(u,\, (dF_2)_n^*(\beta)\bigr)\in (D_N)_n\,\Bigr\}. 
\]
From the assumption that each $F_k:(N, D_N)\to (M_k, D_{M_k})$ is a Dirac map, it follows that 
\[
 \mathrm{pr}_2\Bigl((L_k)_{(n,F_k(n))}\Bigr) = (D_{M_1})_{F_k(n)}\quad (k=1,2). 
\]
This shows that (C) holds. If a point $\bigl(u;\,(dF_k)_n(u),\,\beta\bigr)\in (L_k)_{(n,F_k(n))}$ 
belongs to the space $T_nN\oplus \{\boldsymbol{0}\}\subset T_nN\oplus T_n^*N$, 
we find that $u\in\ker(df_k)_n$ and $\beta=\boldsymbol{0}$. Since the condition 
(\ref{sec3:strong Dirac}), we have 
$(u, \boldsymbol{0})\in \ker(dF_k)_n\cap\ker(D_N)_n=\{\boldsymbol{0}\}$. 
This implies $u=\boldsymbol{0}$. Therefore, (A) holds. 
For any $n\in N$, the each space $\mathrm{pr}_1\bigl((L_k)_{(n,F_k(n))}\bigr)$ coincides with 
$\mathrm{pr}_1\bigl((D_N)_n\cap \bigl(T_nN\oplus \mathrm{Im}\,(dF_k)_n^*\bigr)\bigr)$. 
Consequently, (B) holds. 
\qquad\qquad\qquad\qquad\qquad\qquad\qquad\qquad\qquad\qquad\qquad\qquad\qquad\qquad\qquad\qquad
\qquad\quad $\Box$
\vspace*{0.5cm}

Suppose that two Lie algebroids $A_1\to M_1$ and $A_2\to M_2$ are quasi-equivalent 
by $M_1\overset{J_1}{\leftarrow}X\overset{J_2}{\to}M_2$, and a pair 
$(L_1, L_2)$ of subbundles $L_1\subset {J_1}^!A_1$ and $L_2\subset {J_2}^!A_2^-$. 
Let us choose any smooth section $\alpha\in\varGamma^\infty(A_1)$. 
From the conditions (A) and (C) in Definition \ref{sec3:dfn}, 
there exists a unique element $u\in T_xX$ such that $(u, \alpha_{J_1(x)})\in (L_1)_{(x, J_1(x))}$  
(see Example \ref{sec3:example4}). 
That is, we have a map 
\begin{equation}\label{sec3:eqn5}
\xi_1:\varGamma^\infty(A_1)\ni \alpha\longmapsto -\hat{\alpha}\in \mathfrak{X}(X). 
\end{equation}
as assigning to $\alpha\in \varGamma^\infty(A_1)$ a unique element $-\hat{\alpha}_x\in T_xX$ 
such as $(\hat{\alpha}_x,\,\alpha_{J_1(x)})\in (L_1)_{(x, J_1(x))}$. 
The map (\ref{sec3:eqn5}) defines a left action of $A_1$. 
A right action $\xi_2$ of $A_2$ is defined in the obvious analogous way. 
It follows from (B) that 
\begin{equation*}
 \bigl\{\,\xi_1(\alpha)_x\,|\,\alpha\in \varGamma^\infty(A_1)\,\bigr\}\,= \, 
  T_x\Bigl(J_2^{-1}(J_2(x))\Bigr).
\end{equation*}
Similarly to this case, the right action $\xi_2$ of $A_2$ yields 
\[
 \bigl\{\,\xi_2(\beta)_x\,|\,\beta\in \varGamma^\infty(A_2)\,\bigr\}\,= \, 
  T_x\Bigl(J_1^{-1}(J_1(x))\Bigr).
\] 

\begin{prop}\label{sec4:prop2}
 Two Lie algebroids $A_1\to M_1$ and $A_2\to M_2$ are quasi-equivalent 
 if and only if 
 there exists a smooth manifold $X$ together with surjective submersions 
 $J_k:X\to M_k~(k=1,2)$ such that 
  \begin{enumerate}[\rm\quad(1)]
   \item $A_1$ has a left action $\xi_1$ on $J_1:X\to M_1$ such that 
    \[\ker\,(dJ_2)_x\, =\, 
          \bigl\{\,\xi_1(\alpha)_x\,|\,\alpha\in \varGamma^\infty(A_1)\,\bigr\}
                                                \quad (\forall x\in X);\]
   \item $A_2$ has a right action $\xi_2$ on $J_2:X\to M_2$ such that 
    \[\ker\,(dJ_1)_x\, =\, 
          \bigl\{\,\xi_2(\beta)_x\,|\,\beta\in \varGamma^\infty(A_2)\,\bigr\}
                                                \quad (\forall x\in X).
    \]
  \end{enumerate}
\end{prop}

\noindent({\it Proof})~The necessary condition for $A_1$ and $A_2$ to be strong Morita equivalent 
follows from the above observation. 
Conversely, assume that there exists such a smooth manifold $X$. We define subbundles $L_1$ 
of ${J_1}^!A_1$ and $L_2$ of $({J_2}^!A_2)^-$ as 
\[
L_1 = \coprod_{x\in X}\Bigl\{\,\bigl(\xi_1(\alpha)_x,\,\alpha_{J_1(x)}\bigr)\,\bigm|\, 
\alpha\in \varGamma^\infty(A_1)\,\Bigr\}.
\]
and 
\[
L_2 = \coprod_{x\in X}\Bigl\{\,\bigl(\xi_2(\beta)_x,\,\beta_{J_2(x)}\bigr)\,\bigm|\, 
\beta\in \varGamma^\infty(A_2)\,\Bigr\},
\]
respectively. The condition (C) in Definition \ref{sec3:dfn} hold obviously. 
If we take a zero section $\alpha\equiv \boldsymbol{0}\in\varGamma^\infty(A_1)$, 
then 
$\xi_1(\alpha)_x = \boldsymbol{0}$. 
This shows that (A) holds. 
The condition (B) is verified by the assumptions that 
the images of the action $\xi_1\,(\xi_2)$ are tangent to $J_2$\,(resp. $J_1$)-fibers. 
\qquad\qquad $\Box$
\vspace*{0.5cm}

Basing on the above discussion, we introduce a new binary relation between integrable 
Lie algebroids. 

\begin{dfn}\label{sec3:strong Morita}
 Suppose that both Lie algebroids $A_1\to M_1$ and $A_2\to M_2$ are integrable. 
 They are said to be strong Morita equivalent if they are quasi-equivalent each other, 
 and satisfy the following conditions{\rm :}
  \begin{enumerate}[\quad\rm (A')]
   \item both the left action $\xi_1$ and the right action $\xi_2$ are complete{\rm ;} 
   \item for any smooth section $\alpha\in\varGamma^\infty(A_1)$ and $\beta\in\varGamma^\infty(A_2)$. 
       \[ [\xi_1(\alpha),\,\xi_2(\beta)] = \boldsymbol{0}.\]
  \end{enumerate}
\end{dfn}

It will be shown that strong Morita equivalence is indeed an equivalence relation 
between integrable Lie algebroids. 

\begin{rmk}\label{sec3:flows}
 The second condition in Definition \ref{sec3:strong Morita} indicates that 
 if $\theta^1_t$ and $\theta^2_t$ are the flows of the vector fields 
 $\xi_1(\alpha)$ and $\xi_2(\beta)$, respectively, then it holds that 
 $\theta^1_t\circ \theta^2_s = \theta^2_s\circ\theta^1_t$ 
 for all $t,s$ for which the flows are defined. 
\end{rmk}

\begin{ex}
If two integrable Poisson manifolds $P_1$ and $P_2$ are Morita equivalent, then 
$T^*P_1\to P_1$ and $T^*P_2\to P_2$ are also strong Morita equivalent. 
Indeed, they are quasi-equivalent (see Example \ref{sec3:example11}). 
A left action of $T^*P_1$ on $S\overset{\tau_1}{\to}P_1$ and a right action of 
$(T^*P_2)^-$ on $S\overset{\tau_2}{\to}P_2$ are given like as the action in 
Example \ref{sec3:Poisson map}. The completeness of Poisson maps $\tau_1$ and $\tau_2$ 
implies that both of the actions are complete~(see \cite{CFint04}). Furthermore, it holds that 
\[
 \bigl[\,\Pi_S^\sharp(\tau_1^*df),\, -\Pi_S^\sharp(\tau_2^*dg)\,\bigr] 
 = \Pi_S\bigl(\,\cdot,\,\Pi_S(\tau_1^*df,\,\tau_2^*dg)\,\bigr) = \boldsymbol{0}, 
\]
since fibers of $\tau_1$ and $\tau_2$ are symplectically orthogonal to one another. 
\end{ex}

An $(A_1, A_2)$-bimodule, denoted by $A_1\overset{J_1}{\leftarrow}X\overset{J_2}{\rightarrow}A_2$, 
is a pair of a complete left module $X\overset{J_1}{\to}M_1$ over $A_1$ 
and a complete right module $X\overset{J_2}{\to}M_2$ over $A_2$ 
which makes $A_1$ and $A_2$ be strong Morita equivalent each other 
as in Definition \ref{sec3:strong Morita}. 
Let us consider an $(A_1, A_2)$-bimodule $A_1\overset{J_1}{\leftarrow}X\overset{J_2}{\rightarrow}A_2$ 
and an $(A_2, A_3)$-bimodule $A_2\overset{K_2}{\leftarrow}Y\overset{K_3}{\rightarrow}A_3$.  
Then, the map 
$\widehat{\xi}_1: \varGamma^\infty(A_1)\rightarrow \mathfrak{X}(X\otimes_{A_2}Y)$,  
\[
\varGamma^\infty(A_1)
\ni \alpha\longmapsto \overline{(\,\xi_1\,(\alpha)_x,\, \boldsymbol{0}\,)}
\in T_{\overline{(x,y)}}(X\otimes_{A_2}Y)
\]
and the map $\widehat{\eta}_3: \varGamma^\infty(A_3)\rightarrow \mathfrak{X}(X\otimes_{A_2}Y)$, 
\[
\varGamma^\infty(A_3)
\ni \beta\longmapsto \overline{(\,\boldsymbol{0},\, \eta_3\,(\beta)_y\,)}
\in T_{\overline{(x,y)}}(X\otimes_{A_2}Y) 
\]
induce a complete left action of $A_1$ on $\widehat{J}_1:X\otimes_{A_2}Y\to M_1,\,
\overline{(x, y)}\mapsto J_1(x)$ and a complete right action 
of $A_3$ on $\widehat{K}_3:X\otimes_{A_2}Y\to M_3,\,\overline{(x, y)}\mapsto K_3(y)$, 
respectively. It is easily verified that those actions satisfy 
\[\ker\,(d\widehat{K}_3)_x\, \supset\, 
   \bigl\{\,\widehat{\xi}_1(\alpha)_x\,|\,\alpha\in \varGamma^\infty(A_1)\,\bigr\}
\]
and 
\[\ker\,(d\widehat{J}_1)_x\, \supset\, 
  \bigl\{\,\widehat{\eta}_3(\beta)_x\,|\,\beta\in \varGamma^\infty(A_3)\,\bigr\}. 
\]
If $\overline{(u,\, v)}$ is any point in $\ker\,(d\widehat{K}_3)_{\overline{(x,y)}}$, then 
there exists a smooth section $\beta\in\varGamma^\infty(A_2)$ such that $v = \eta_2(\beta)_y$. 
Consequently, we have 
\begin{align*}
 (dJ_2)_x\bigl(u - \xi_2(\beta)_y\bigr) &= (dJ_2)_x(u) - \rho_2(\beta) = (dK_2)_y(v) - \rho_2(\beta)\\
  &= (dK_2)_y(\eta_2(\beta)_y) - \rho_2(\beta) = \boldsymbol{0}. 
\end{align*}
That is, $u-\xi_2(\beta)\in\ker\,(dJ_2)_x$. By the assumption, there exists a smooth section 
$\alpha\in\varGamma^\infty(A_1)$ such that $\xi_1(\alpha)_x = u - \xi_2(\beta)_x$. 
Therefore, 
\[
 (u,\,v) = (u-\xi_2(\beta)_x,\, \boldsymbol{0}) + (\xi_2(\beta)_x,\,\eta_2(\beta)_y).
\]
This implies that $\overline{(u,\,v)} = \overline{(\xi_1(\alpha),\,\boldsymbol{0})}$. 
As a result, we show that 
\[\ker\,(d\widehat{K}_3)_x\, =\, 
   \bigl\{\,\widehat{\xi}_1(\alpha)_x\,|\,\alpha\in \varGamma^\infty(A_1)\,\bigr\}.
\]
Similarly, 
\[\ker\,(d\widehat{J}_1)_x\, =\, 
  \bigl\{\,\widehat{\eta}_3(\beta)_x\,|\,\beta\in \varGamma^\infty(A_3)\,\bigr\}. 
\]
The observation leads us to the conclusion 
that the leaf space $X\otimes_{A_2}Y$ is an $(A_1, A_3)$-bimodule. 

\begin{ex}\label{sec3:example14}
If $X\to M$ is the right module over $A$, 
then $\{\ast\}\leftarrow X\to M$ is the $(\ast, A)$-bimodule. 
Similarly, $M\leftarrow X\to \{\ast\}$ turns out to be the $(A, \ast)$-bimodule 
if $X\to M$ is the left module over $A$. 
\end{ex}

\noindent On the basis of those observations, we can show the following proposition. 

\begin{prop}\label{sec3:equivalence relation}
 Strong Morita equivalence for integrable Lie algebroids is an equivalence relation. 
\end{prop}

\noindent({\it Proof})~
The transitivity holds obviously by the above observation. 
Let $A\to M$ be an integrable Lie algebroid and 
$\varGamma(A)\rightrightarrows M$ a Lie groupoid integrating $A\to M$. 
From Example \ref{sec3:groupoid action}, we have the left action $\xi$ by 
(\ref{sec3:groupoid action1}) and the right action $\eta$ by (\ref{sec3:groupoid action2}). 
It is obvious that those actions are complete. 
As for the left action $\xi$, we have 
\[
 (d\boldsymbol{s})_{\varepsilon(x)}\bigl(\xi(\alpha)\bigr) 
   = \left.\frac{d}{dt}\right|_{t=0}\boldsymbol{s}(\gamma(t)\cdot h)
   = \left.\frac{d}{dt}\right|_{t=0}\boldsymbol{s}(h) = \boldsymbol{0}\quad (\forall x\in M). 
\]
Similarly, 
\[
 (d\boldsymbol{t})_{\varepsilon(x)}\bigl(\eta(\beta)\bigr) 
   = \left.\frac{d}{dt}\right|_{t=0}\boldsymbol{t}(g\cdot\delta(t))
   = \left.\frac{d}{dt}\right|_{t=0}\boldsymbol{t}(h) = \boldsymbol{0}\quad (\forall x\in M). 
\]
From this, it follows that $\ker(d\boldsymbol{s})_{\varepsilon(x)}
=\{\,\xi(\alpha)_x\,|\,\alpha\in\varGamma^\infty(A)\,\}$ and $\ker(d\boldsymbol{t})_{\varepsilon(x)} 
=\{\,\eta(\beta)_x\,|\,\beta\in\varGamma^\infty(A)\,\}$. 
Moreover, 
\begin{align*}
 (d\boldsymbol{s})_{\varepsilon(x)}\bigl(\eta(\beta)\bigr) 
   &= \left.\frac{d}{dt}\right|_{t=0}\boldsymbol{s}(g\cdot\delta(t))
    = \left.\frac{d}{dt}\right|_{t=0}\boldsymbol{s}(\delta(t))\\
   &= (d\boldsymbol{s})_{\varepsilon(x)}(\widehat{\beta}_{\varepsilon(x)}) 
    = (d\boldsymbol{s})_{\varepsilon(x)}
   \left((dR_{\varepsilon(x)})(\beta_{\varepsilon(x)})\right)
    = \boldsymbol{0},  
\end{align*}
since the right invariant vectors $\widehat{\beta}$ lie in the $\boldsymbol{s}$-fibers. 
Therefore, we have $[\xi(\alpha),\, \eta(\beta)]=\boldsymbol{0}$ for any 
$\alpha, \beta\in\varGamma^\infty(A)$. This results in that an integrable Lie algebroid 
$A$ is strong Morita equivalent to itself. 
Lastly, 
suppose that $A_1$ and $A_2$ are strong Morita equivalent by $(A_1, A_2)$-bimodule 
$A_1\overset{J_1}{\leftarrow}X\overset{J_2}{\rightarrow}A_2$. 
Defining a left action $\xi'$ of $A_2$ and a right action $\eta'$ of $A_1$ as 
$\xi'(\beta):= -\xi_2(\beta)~(\forall\beta\in\varGamma^\infty({J_2}^*A_2))$ 
and $\eta' (\alpha):= -\xi_1(\alpha)~(\forall\alpha\in\varGamma^\infty({J_1}^*A_1))$, 
respectively, we obtain an $(A_2, A_1)$-bimodule 
$A_2\overset{J_2}{\leftarrow}X\overset{J_1}{\rightarrow}A_1$. 
This shows that the reflectivity holds. 
\qquad\qquad\qquad\qquad\qquad\qquad\qquad\qquad\qquad\qquad\qquad\quad $\Box$

%%%%%%%%%%%%%%%%%%%%%%%%%%%%%%%%%%%%%%%%%%%%%%%%%%%%%%%%%%%%%%%%%%%%%%%%%%%%%
\section{Equivalence of the categories of infinitesimal actions}
\subsection{$A$-paths}
%%%%%%%%%%%%%%%%%%%%%%%%%%%%%%%%%%%%%%%%%%%%%%%%%%%%%%%%%%%%%%%%%%%%%%%%%%%%

Before proceeding to the main theorem, let us review briefly the basics of an $A$-path. 
For further details, we refer to M. Crainic and R. J. Fernandes~\cite{CFint03}, \cite{DZpoi05}. 

\begin{dfn}
 Let $A\overset{\pi}{\to} M$ be a Lie algebroid with an anchor map $\rho:A\to TM$. 
 An $A$-path is a smooth path $a:I\to A$ which projects 
 to a base path $\pi\circ a:I\to M$ such that 
 \[
  \rho\bigl(a(t)\bigr) \,=\, \frac{d}{dt}\pi\bigl(a(t)\bigr)\quad (\forall t\in I). 
 \]
 Here, $I=[0,\,1]$ stands for the unit interval. 
\end{dfn}

A map $a_\epsilon(t):=a(t,\epsilon):I\times I\to A$ is called a variation of $A$-paths if 
$a_\epsilon$ is a family of $A$-paths of class $C^2$ on $\epsilon$ with the property that 
the base paths $c_\epsilon(t)=\pi\circ a_\epsilon(t):I\times I\to M$ have fixed end points. 
For a time-dependent section $\sigma_\epsilon$ of $A$ such that 
$\sigma_\epsilon(t,\,c_\epsilon(t)) = a_\epsilon(t)$, consider a function 
\[
 b(\epsilon, t) = \int_0^t\phi_{\sigma_\epsilon}^{t,s}\frac{d\sigma_\epsilon}{d\epsilon}
 (s,\,c_\epsilon(s))\,ds, 
\]
where $\phi_{\sigma_\epsilon}^{t,s}$ denotes the flow of the time-depending section 
$\sigma_\epsilon$. 

\begin{dfn}
 Two $A$-path $a_0$ and $a_1$ are $A$-homotopic if there exists a variation $a_\epsilon$ 
 such that $b(\epsilon, 1)=0$ for $\varepsilon\in I$. We write $a_0\sim a_1$ when 
 $a_0$ and $a_1$ are $A_1$-homotopic. 
\end{dfn}

Denote by $P(A)$ the space of $A$-paths of class $C^1$ for a Lie algebroid $A$. 
Suppose that the Lie algebroid $A\to M$ is integrable. Then, the quotient 
$\mathcal{G}(A) := P(A)/\sim$ is a smooth manifold~(see \cite{CFint03}). 
As a matter of fact, $\mathcal{G}(A)$ turns out to be a Lie groupoid 
over $M$ with the structure maps $\boldsymbol{s}([a]):=\pi(a(0))$ and 
$\boldsymbol{t}([a]):=\pi(a(1))$~(see \cite{CFint03}). 
This Lie groupoid is called a Weinstein groupoid.

%%%%%%%%%%%%%%%%%%%%%%%%%%%%%%%%%%%%%%%%%%%%%%%%%%%%%%%%%%%%%%%%%%%%%%%%%
\subsection{Main result}
%%%%%%%%%%%%%%%%%%%%%%%%%%%%%%%%%%%%%%%%%%%%%%%%%%%%%%%%%%%%%%%%%%%%%%%%%

In algebra, Morita equivalence implies an equivalence of the 
categories of modules. To be concrete, if two algebras $R,\,S$ are Morita equivalent, 
then $\mathcal{M}_R$ and $\mathcal{M}_S$ are equivalent. Here, $\mathcal{M}_R$ denotes 
the category of right $R$-modules, whose objects are right $R$-modules and whose morphisms 
between $M_1$ and $M_2$ are $R$-homomorphisms. 
Basing on this well-known fact, we introduce the category of modules over Lie algebroids 
as follows: 

\begin{dfn}
Let $(A\to M,\,\llbracket\cdot, \cdot\rrbracket,\,\rho)$ be an integrable Lie algebroid. 
The category of modules over a Lie algebroid is the category $\mathcal{M}\,(A)$ 
whose objects are right modules over $A$ 
and whose morphisms between $\mu:N\to M$ and $\mu':N'\to M$ are smooth map 
$f:N\to N'$ such that $\mu'\circ f=\mu$ and, for each $\alpha\in \varGamma^\infty(A)$, 
their respective vector fields $\xi(\alpha)\in\mathfrak{X}(N)$ and $\xi'(\alpha)\in\mathfrak{X}(N')$ 
are $f$-related\rm{:}
\[
\xi'(\alpha)_{f(n)} = (df)_n(\xi(\alpha)_n)\quad (\forall n\in N) 
\]
\end{dfn}

Let $(A_1\to M_1,\,\llbracket\cdot, \cdot\rrbracket_1,\, \rho_1)$ 
and $(A_2\to M_2,\,\llbracket\cdot, \cdot\rrbracket_2,\, \rho_2)$ be integrable Lie algebroids and 
assume that $A_1$ and $A_2$ are strong Morita equivalent. 
Let $N\overset{\mu}{\to} M_1$ be a right module over $A_1$. From the assumption, 
there exists an $(A_1, A_2)$-bimodule 
$A_1\overset{J_1}{\leftarrow}X\overset{J_2}{\rightarrow}A_2$. We remark that 
$A_1$ acts on $\mu:N\to M_1$ and $J_1:X\to M_1$ from the right and the left, respectively. 
The right module $N\to M_1$ over $A_1$ is thought of a $(\ast, A_1)$-bimodule 
(see Example \ref{sec3:example14}). Hence, as discussed earlier, the tensor product 
\begin{equation*}
\nu: N\otimes_{A_1}X\longrightarrow M_2,\quad 
\overline{(\,n,\, x\,)}\longmapsto J_2(x)
\end{equation*} 
turns out to be a $(\ast, A_2)$-module, that is, a right module over $A_2$. 
In addition, given a morphism $f:N\to N'$ in $\mathcal{M}\,(A_1)$, we define the map 
$\widehat{f}:N\otimes_{A_1}X\to N'\otimes_{A_1}X$ as 
$\widehat{f}\,(\overline{n, x}) = (\overline{f(n), x})$ 
for any $(\overline{n, x})\in N\otimes_{A_1}X$. 
The map $\widehat{f}$ turns out easily to satisfy $\nu'\circ\widehat{f}=\nu$. 
As a result, we obtain a functor $\mathcal{S}$ from $\mathcal{M}\,(A_1)$ to 
$\mathcal{M}\,(A_2)$ which assigns to each object $N\to M_1$ 
an object $N\otimes_{A_1}X\to M_2$, 
and to each morphism $f:N\to N'$ in $\mathcal{M}\,(A_1)$ a morphism 
$N\otimes_{A_1}X\to N'\otimes_{A_1}X$. 
\vspace*{0.5cm}

\noindent
This observation leads us to the following proposition. 

\begin{prop}\label{sec4:covariant functor}
If $A_1$ and $A_2$ are strongly Morita equivalent, then 
there exists a covariant functor from $\mathcal{M}\,(A_1)$ to $\mathcal{M}\,(A_2)$. 
\end{prop}

\noindent 
In a similar way, we can obtain a covariant functor $\mathcal{T}$ from $\mathcal{M}\,(A_2)$ to 
$\mathcal{M}\,(A_1)$. 
\vspace*{0.5cm}

Suppose that $A_1$ and $A_2$ be strong Morita equivalent by an $(A_1, A_2)$-bimodule 
$A_1\overset{J_1}{\leftarrow}X\overset{J_2}{\rightarrow}A_2$. 
According to the observation in the previous section, one can obtain a $(A_1, A_1)$-bimodule 
$A_1\overset{\widehat{J_1}}{\leftarrow}X\otimes_{A_2}X\overset{\widehat{J}_1}{\rightarrow}A_1$. 
Denote by $\zeta$ the right action of $A_1$ on $J_1:X\to M_1$ and 
let $t\mapsto a(t)$ be an $A_1$-path with a base path $t\mapsto c(t):=\pi(a(t))$ 
starting at $m\in M$. 
For a point $x\in J_1^{-1}(m)\subset X$, 
we consider the following 
differential equation with the initial value problem: 
\begin{equation}\label{sec4:ode}
 \frac{d}{dt}u(t) = \zeta_{u(t)}\bigl(a(t)\bigr),\quad u(0)=x. 
\end{equation}
To verify that the equation (\ref{sec4:ode}) has a unique solution defined on the entire  
unit interval, we choose a time-dependent smooth section $\sigma$ of $A_1$ 
which has compact support, and satisfies $\sigma\bigl(t, c(t)\bigr)=a(t)$. 
A solution of (\ref{sec4:ode}) is an integral curve of a time-dependent 
vector field 
\[
V_{(t, x)}:=\zeta_x\bigl(\sigma(t, J_1(x))\bigr)\quad \bigl(\,(t,x)\in I\times X\,\bigr) 
\] 
induced from $\sigma$.  
Conversely, suppose that $u$ is an integral curve of $V$. 
Then, it is verified that 
$J_1\circ u$ is an integral curve of $\rho_1(\sigma)$ with the initial point $m$. 
Indeed, we have 
\[
 \frac{d}{dt}(J_1\circ u)(t) 
 = (dJ_1)_{u(t)}\Bigl(\zeta_{u(t)}\bigl(\sigma\bigl(t,\, J_1\circ u(t)\bigr)\bigr)\Bigr) 
 = \rho_1\bigl(\sigma(t,\, J_1\circ u(t))\bigr). 
\]
by (\ref{sec3:eqn1}) and $J_1(u(0))=J_1(x)=m$. 
Consequently, the curve $J_1\circ u$ coincides with the base path $c$. 
From this, it follows that 
\begin{equation}
\frac{d}{dt}u(t) = \zeta_{u(t)}\bigl(\sigma\bigl(t,\, c(t)\bigr)\bigr) 
= \zeta_{u(t)}\bigl(a(t)\bigr). 
\end{equation}
Therefore, the equation (\ref{sec4:ode}) has a unique solution. 
Moreover, $u$ is defined on the entire $I$ since the completeness of the action $\zeta$ 
implies that $\zeta(\sigma(t,\,\cdot))$ is complete whenever 
$\sigma(t,\,\cdot)$ has compact supported. 
Now, let us take the homotopy class $[a]\in\mathcal{G}(A_1)$ for any $A_1$-path $a$ 
and any point $x\in X$ such that $J_1(x)=\pi(a(0))$. 
Then, the map 
\begin{equation}\label{sec4:integral curve}
 \mathcal{G}(A_1)\ni [a] \longmapsto \overline{(x', u(1))}\in X\otimes_{A_2}X, 
\end{equation} 
which assigns to $[a]$ an equivalent class of a point $(x', u(1))$ along the integral curve 
$u$ determined by the above observation, is a diffeomorphism between 
$\mathcal{G}(A_1)$ and $X\otimes_{A_2}X$. It can be shown that the map is well-defined by 
using the method similar to the proof of Lemma 2 in \cite{CFint04}. 

\begin{rmk}
We remark that a point $(x, u(1))$ belongs to the fiber-product $X\times _{M_2}X$. 
Indeed, it follows from (\ref{sec4:ode}) and the assumption that 
\[
 \left.\frac{d}{dt}\right|_t J_2(u(t)) = (dJ_2)_{u(t)}\bigl(\zeta(a(t))\bigr) = \boldsymbol{0}. 
\]
In other words, $J_2(u(t))=\text{const.}$ for each $t\in I$. 
Therefore, $J_2(u(1))=J_2(u(0))=J_2(x)$. 
\end{rmk}

\begin{thm}\label{sec4:main result}
 If $A_1$ and $A_2$ are strongly Morita equivalent, then their categories of modules 
 $\mathcal{M}\,(A_1)$ and $\mathcal{M}\,(A_2)$ are equivalent. 
\end{thm}

\noindent({\it Proof})~ 
Let $N\overset{\mu}{\to} M_1$ and $N'\overset{\mu'}{\to} M_1$ be objects in $\mathcal{M}\,(A_1)$, 
and $f:N\to N'$ a morphism in $\mathcal{M}\,(A_1)$. 
From Proposition \ref{sec4:covariant functor}, there exist 
the covariant functors $\mathcal{S}:\mathcal{M}\,(A_1)\to \mathcal{M}\,(A_2)$ and 
$\mathcal{T}:\mathcal{M}\,(A_2)\to \mathcal{M}\,(A_1)$, and then, one can obtain 
two right modules 
\[
\mathcal{T}\circ\mathcal{S}\,(N)=(N\otimes_{A_1}X)\otimes_{A_2}X \longrightarrow M_1, 
\quad 
(\overline{(\overline{n, x'})_1,\,x})_2\longmapsto J_1(x), 
\] 
and 
\[
\mathcal{T}\circ\mathcal{S}\,(N')=(N'\otimes_{A_1}X)\otimes_{A_2}X \longrightarrow M_1, 
\quad 
(\overline{(\overline{n', x'})_1,\,x})_2\longmapsto J_1(x), 
\]
over $A_1$, and a morphism
\[
 \mathcal{T}\circ \mathcal{S}\,(f): \mathcal{T}\circ \mathcal{S}\,(N)
 \to \mathcal{T}\circ \mathcal{S}\,(N'), 
 \quad (\overline{(\overline{n', x'})_1,\,x})_2\mapsto (\overline{(\overline{f(n), x'})_1,\,x})_2
\] 
from $\mathcal{T}\circ \mathcal{S}~(N)$ to $\mathcal{T}\circ \mathcal{S}~(N')$. 
Here, $A_1\overset{J_1}{\leftarrow}X\overset{J_2}{\rightarrow}A_2$ is an $(A_1, A_2)$-bimodule. 

On the other hand, for any right module $N\overset{\mu}{\to} M_1$ over $A_1$, 
we define a map $\Psi_N:\mathcal{T}\circ \mathcal{S}\,(N)\to N$ as 
\[
\Psi_N: (N\otimes_{A_1}X)\otimes_{A_2}X\overset{\simeq}{\longrightarrow} N,\quad 
(\overline{(\overline{n, x'})_1,\,x})_2\longmapsto (\overline{x', x})_2\cdot n, 
\]
where the element $(\overline{x', x})_2$ in 
$X\otimes_{A_2}X$ is thought of an element in $\mathcal{G}(A_1)$ 
by (\ref{sec4:integral curve}), and where $(\overline{x, x'})_2\cdot n$ means the point $\phi_1(n)$ 
on the integral curve $\phi_t(n)$ starting at $n\in N$ which is determined by 
$[b]:=(\overline{x, x'})_2\in\mathcal{G}(A_1)$. Remark that the map $\Psi_N$ is well-defined 
by Remark \ref{sec3:flows}. 
Since their respective vector fields induced those actions $\varrho$ and $\varrho'$ of $A_1$ 
are $f$-related, we have 
\[
 \frac{d}{dt}(f\circ \phi_t)(n) = (df)_{\phi_t(n)}\left(\frac{d}{dt}\phi_t(n)\right) 
 = (df)_{\phi_t(n)}\bigl({\varrho}(b(t))_{b(t)}\bigr) = \varrho'_{f(b(t))}(b(t)). 
\]
That is, $t\mapsto f\circ \phi_t(n)$ is an integral curve of $\varrho'(b(t))$. 
From the uniqueness, it follows that $f((\overline{x', x})_2\cdot n)=(\overline{x', x})_2\cdot f(n)$. 
In other words, the diagram 
\[\begin{CD}
\mathcal{T}\circ\mathcal{S}\,(N)  @>{\Psi_N}>> N\\
@V{\mathcal{T}\circ\mathcal{S}\,(f)}VV           @VV{f}V \\
\mathcal{T}\circ\mathcal{S}\,(N') @>>{\Psi_{N'}}> N'.
\end{CD}
\]
commutes. Consequently, the functor $\mathcal{T}\circ\mathcal{S}$ is natural isomorphic to 
the identity functor $\mathrm{Id}_{\mathcal{M}(A_1)}$. 
Similarly, it can be shown that there exists also a natural isomorphism 
between $\mathcal{S}\circ\mathcal{T}$ and $\mathrm{Id}_{\mathcal{M}(A_2)}$. 
This completes the proof. 
\qquad\qquad\qquad\qquad\qquad\qquad\qquad\qquad\qquad\qquad\qquad 
\qquad\qquad\qquad$\Box$ 

%%%%%%%%%%%%%%%%%%%%%%%%%%%%%%%%%%%%%%%%%%%%%%%%%%%%%%%
\subsection{Application to the Hamiltonian category}
%%%%%%%%%%%%%%%%%%%%%%%%%%%%%%%%%%%%%%%%%%%%%%%%%%%%%%%
Let $(M, D_M)$ be a Dirac manifold. 
The Hamiltonian category \cite{BCdir08} of $(M, D_M)$, 
denoted by $\overline{\mathcal{M}}\,(M, D_M)$, 
is a category whose objects are strong Dirac maps 
$F:N\to M$ and whose morphisms are forward Dirac maps $\varphi:N\to N'$ satisfying 
$F=F'\circ\varphi$. 

Let us focus on the case where $(M, D_M)$ is integrable. 
We define the subcategory $\mathcal{H}\,(M, D_M)$ of 
$\overline{\mathcal{M}}\,(M, D_M)$ as the category whose objects are complete strong Dirac maps 
$F:N\to M$ which are surjective submersions and whose morphisms are forward Dirac maps 
$\varphi:N\to N'$ such that $F=F'\circ\varphi$. 
Suppose that we are given a morphism $\varphi:N\to N'$ in $\mathcal{H}\,(M, D_M)$. 
As mentioned in Section 3, 
each strong Dirac map $N\overset{F}{\to}M$ and $N'\overset{F'}{\to}M$ 
induces an infinitesimal action $\zeta:\varGamma^\infty(D_M)\to \mathfrak{X}(N)$ 
and $\zeta':\varGamma^\infty(D_M)\to \mathfrak{X}(N')$ by (\ref{sec3:Dirac action}).  
That is, for $(Y, \beta)\in \varGamma^\infty(D_M)$, we can obtain 
$X_n\in T_nN$ and $X'_{n'}\in T_{n'}N'$ 
such that $Y_{F(n)}=(dF)_n(X_n)$ and 
$Y_{F'(n')}=(dF')_{n'}(X_{n'})$. Note that 
$Y_{F(n)} = d(F'\circ\varphi)_n(X_n) = (dF')_{\varphi(n)}\bigl((d\varphi)_n(X_n)\bigr)$ 
by $F=F'\circ\varphi$ and we have 
$(dF')_{\varphi(n)}\bigl(X'_{\varphi(n)}-(d\varphi)_n(X_n)\bigr)=\boldsymbol{0}$. 
That is, $X'_{\varphi(n)}-(d\varphi)_n(X_n)\in\ker\,(dF')_{\varphi(n)}$. 
In addition, $X'_{\varphi(n)}-(d\varphi)_n(X_n)$ is also an element in 
$\ker\,(D_{N'})$ since $\varphi$ is a forward-Dirac map. Consequently,  
\[
X'_{\varphi(n)}-(d\varphi)_n(X_n)\in \ker\,(dF')_{\varphi(n)}\cap\ker\,(D_{N'})_{\varphi(n)}
=\{\boldsymbol{0}\}. 
\] 
This shows that the vector fields $\zeta(Y,\beta)$ and $\zeta'(Y,\beta)$ are $\varphi$-related. 
Therefore, the subcategory $\mathcal{H}\,(M, D_M)$ of the Hamiltonian category 
$\overline{\mathcal{M}}\,(M, D_M)$ can be regarded as the category of modules 
$\mathcal{M}\,(D_M)$ over $D_M$. 
\vspace*{0.5cm}

For a closed 2-form $B$ on $M$, a subbundle 
\[
 \tau_B(D_M) := \bigl\{\,(Y,\, \beta + \mathrm{i}_YB)\,|\,(Y,\,\beta)\in D_M\,\bigr\}
 \subset TM\oplus T^*M. 
\]
satisfies the conditions in Example \ref{sec2:example of Lie algebroid}. 
In other words, $\tau_B(D_M)$ is a Dirac structure on $M$. 
The Dirac structure $\tau_B(D_M)$ associated to a closed 2-form $B$ on $M$ is called a 
gauge transformation by $B$~(see H. Bursztyn and O. Radko \cite{BRgau03}). 
In Example 6.6~\cite{BCWZint04}, it is proven that $\tau_B(D_M)$ is integrated to 
the Lie groupoid $\mathcal{G}(D_M)\rightrightarrows M$ associating with $D_M$. 
By noting Example \ref{sec3:groupoid action} and Proposition \ref{sec3:equivalence relation}, 
it is easily show the following proposition: 

\begin{prop}
 If $(M,\,D_M)$ is an integrable Dirac manifold, then two Lie algebroids $D_M$ and 
 $\tau_B(D_M)$ are strong Morita equivalent each other. 
\end{prop}

By using this proposition and Theorem \ref{sec4:main result}, 
we can recover partially Proposition 2.8 \cite{BCdir08}. 

\begin{cor}{\rm (cf. H. Bursztyn and M. Crainic \cite{BCdir08})}~
Let $D_M$ be an integrable Dirac structure and $B$ a closed 2-form on $M$. 
Then, the subcategories $\mathcal{H}\,(M. D_M)$ and $\mathcal{H}\,(M. \tau_B(D_M))$ 
are equivalent. 
\end{cor}

%%%%%%%%%%%%%%%%%%%%%%%%%%%%%%%%%%%%%%%%%%%%%%%%%%%%%%%


\begin{thebibliography}{99}
\bibitem{BCdir05} H. Bursztyn and M. Crainic: {\it Dirac structures, moment maps and quasi-Poisson 
manifolds}. In: {\it The breadth of symplectic and Poisson geometry. Festschrift in honor of 
Alan Weinstein}, 1--40, Progr. Math {\bf 232}, Birkh$\ddot{\rm a}$user Boston, 2005. 

\bibitem{BCdir08} H. Bursztyn and M. Crainic: {\em Dirac geometry, quasi-Poisson actions and 
$D/G$-valued moment maps}, J. Differential Geom. {\bf 82}, 2009, no. 3, 501--566. 

\bibitem{BCWZint04} H. Bursztyn, M. Crainic, A. Weinstein, and C. Zhu: 
{\em Integration of twisted Dirac brackets}, Duke Math. J. {\bf 123}, 2004, no. 3, 549--607. 

\bibitem{BRgau03} H. Bursztyn and O. Radko: 
{\em Gauge equivalence of Dirac structures and symplectic groupoids}, 
Ann. Inst. Fourier (Grenoble) {\bf 53}, 2003, 309--337. 

\bibitem{BWpoi05} H. Bursztyn and A. Weinstein: 
{\it Poisson geometry and Morita equivalence}. In: {\it Poisson Geometry, 
Deformation Quantisation and Group Representations}, 1--78, 
London Math. Soc. Lecture Notes Series {\bf 323}, Cambridge University Press, Cambridge, 2005. 

\bibitem{CFint03} M. Crainic and R. L. Fernandes: {\it Integrability of Lie brackets}, 
Ann. of Math. (2) {\bf 157}, 2003, no. 2, 575--620. 

\bibitem{CFint04} M. Crainic and R. L. Fernandes: {\it Integrability of Poisson brackets}, 
J. Differential Geom. {\bf 66}, 2004, no. 1, 71--137. 

\bibitem{DZpoi05} J.-P. Dufour and N. T. Zung: {\it Poisson Structures and Their Normal Forms}. 
Prog. in Math. {242}, Birkh$\ddot{a}$user, Bassel, 2005.

\bibitem{HMalg90} P. Higgins and K. Mackenzie: {\it Algebraic constructions in the category of 
Lie algebroids}, Journal of Algebra {\bf 129}, 1990, 194--230. 

\bibitem{Lqua01} N. P. Landsman: {\it Quantized reduction as a tensor product}. 
In: {\it Quantization of Singular Symplectic Quotients}, N. P. Landsman, M. Pffaum, 
M. Schlichenmaier, eds., 137--180, Birkh$\ddot{\rm a}$user, Bassel, 2001.

%\bibitem{LSXnon09} C. Laurent-Gengoux, M. Sti$\acute{\rm e}$non and P. Xu: 

\bibitem{Mgen05} K. Mackenzie: {\it General theory of Lie groupoids and Lie algebroids}. 
London Math. Soc. Lecture Notes Series {\bf 213}, Cambridge University Press, Cambridge, 2005. 

\bibitem{Mdua58} K. Morita: {\it Duality for modules and its applications to the theory 
of rings with minimum condition}. Sci. Rep. Tokyo Kyoiku Daigaku Sect. A {\bf 6}, 1958, 83--142. 

\bibitem{Wgeo05} A. Weinstein: {\it The geometry of momentum}, 
G$\acute{\rm e}$om$\acute{\rm e}$trie au XX$\grave{\rm e}$me Si$\grave{\rm e}$cle, Histoire et Horizons, 
Hermann, Paris, 2005. Arxiv: Math.SG/0208108. 

\bibitem{Xmor91} P. Xu: {\it Morita equivalence of Poisson manifolds}. 
Comm. Math. Phys. {\bf 142}, 1991, 493--509. 


\end{thebibliography}
\end{document}